\documentclass[12pt]{article}
\usepackage{amssymb}
\def\qed{\hfill$\Box$\par}
\def\a{\alpha}
\def\rla{\leftrightarrow}

\def\l{\lambda}

\def\dis{\displaystyle}
\def\cl{\centerline}

\def\ul{\underline}
\def\wt{\widetilde}

\def\bs{\backslash}
\def\hs{\hspace*}
\def\vs{\vspace*}
\def\rb{\raisebox}
\def\ra{\rangle}

\def\la{\langle}
\def\ni{\noindent}
\def\ptl{\partial}

\def\hi{\hangindent}
\def\ha{\hangafter}

\def\Z{\mathbb{Z}}
\def\Q{\mathbb{Q}}

\def\C{\mathbb{C}}

\def\jj{{\ul j}}
\def\ii{{\ul i}}\def\nn{{\ul n}}\def\mm{{\ul m}}\def\kk{{\ul k}}%
\begin{document}
%
%\par\ Must put ZHU as the corresponding author
%\setcounter{page}{0}\pagebreak\par
%
\def\LIE{\mbox{${\rm Der\,}\C_q$}}
\cl{\large\bf Representations of the Derivation Algebra of the}
\cl{\large\bf Localization of the Quantum Plane at $q=-1$} \vskip5pt
\cl{(to appear in {\it Comm. Alg.})}

\vs{10pt}\par
\cl{Yucai Su$^{\,1}$
 and Linsheng Zhu$^{\,2}$}\par
\cl {\small $^{1}\,$Department of Mathematics,
                Shanghai Jiaotong University, Shanghai 200030,
                P.~R.~China}
\cl{\small Email: ycsu@sjtu.edu.cn}\par
\cl{\small $^{2}\,$Department of Mathematics,
 Changshu Institute of Technology, Jiangsu 215500, P.~R.~China}
\vskip12pt%
\begin{minipage}[rm]{5.8truein}
\parskip .05 truein \baselineskip 4pt \lineskip 4pt
\par\ \vs{-5pt}\par\ni
{\small %
\cl{\bf ABSTRACT}
\par\ni
We determine the irreducible weight modules with
weight multiplicities at most 1 over the derivation algebra
of the localization of the quantum plane at $q=-1$.
\vs{8pt}\par\ni%
\vs{8pt}\par\ni%
}
\end{minipage}
\par\
\vs{5pt}\par \cl {\bf
             1. \ INTRODUCTION
} \par%\ni
Defined as the noncommutative associative algebra
$\C_q=\C_q[x_1,x_2,x_1^{-1},x_2^{-1}]$ of (noncommutative) Laurent
polynomials over the complex field $\C$ with the relation
$x_2x_1=qx_1x_2$ (where $q\ne0,1$ is a fixed complex number),
the {\it quantum plane},
which has drawn some
\mbox{authors'} attentions (eg, [1,\,2,\,4,\,9,\,10]), is a
\mbox{fundamental} ingredient in quantum groups (cf.~[3]). Under
the usual commutator, it is a Lie algebra, denoted as $\C_q^L$,
with basis
$$B_0=\{x^{\mm}=x_1^{m_1}x_2^{m_2}\,|\,\mm=(m_1,m_2)\in\Z^2\},$$ and
relation
$$
[x^{\mm},x^{\nn}]=(q^{m_2n_1}-q^{m_1n_2})x^{\mm+\nn}\mbox{ \ \ for
\ \ }\mm,\nn\in\Z^2.\eqno(1.1)
$$
It is the $q$-analog of the Lie algebra with basis $B_0$ and
relation
$$
[x^{\mm},x^{\nn}]=(m_2n_1-m_1n_2)x^{\mm+\nn}\mbox{ \ \ for \ \
}\mm,\nn\in\Z^2.\eqno(1.2)
$$
The Lie algebra with relation (1.2) is usually referred to as the
{\it Virasoro-like algebra} in the literature (cf.~[2,\,9]). Thus
the Lie algebra $\C^L_q$ is also referred to as the {\it
$q$-analog of the Virasoro-like algebra} [2,\,4,\,9]. One has
$$
\C^L_q=\C'_q\oplus Z(\C_q),\eqno(1.3)
$$%
where%
$$\C'_q=[\C^L_q,\C^L_q],\,\ \ Z(\C_q)={\rm
span}\{x^{\mm}\,\,|\,\,\mm\in\Z^2,q^{m_1}=q^{m_2}=1\},$$%
are respectively the {\it derived subalgebra} of $\C^L_q$ and the
{\it center of $\,\C_q$}. Obviously, $\C'_q$ is a simple Lie
subalgebra of $\C_q$. It is straightforward to compute (see also,
[1,\,2]) that
$$
\LIE={\rm ad\,}\C'_q\oplus (Z(\C_q)\ptl_1\oplus Z(\C_q)\ptl_2),\ \
\LIE\cong{\rm Der\,}\C'_q,\eqno(1.4)
$$
where $\LIE$ is the {\it derivation algebra of the associative
algebra $\C_q$}, ${\rm Der\,}\C'_q$ is the {\it derivation algebra
of the Lie algebra $\C'_q$}, ${\rm ad\,}\C'_q=\{{\rm
ad\,}u\,|\,u\in\C'_q\}$ (which is in fact isomorphic to $\C'_q$)
is the {\it inner derivation algebra of $\,\C'_q$}, and where
$\ptl_1,\ptl_2$ are {\it degree derivations of $\,\C_q$} defined
by
$$%
\ptl_s(x^{\mm})=m_sx^{\ul m}\mbox{ \ \ for \ \ }\mm\in\Z^2,\ \
s=1,2.\eqno(1.5)
$$%
The isomorphism in (1.4) is given by the restriction $d\mapsto
d|_{\C'_q}$ for $d\in\LIE$. By (1.3), one sees that $\C'_q$ is a
centerless perfect Lie algebra (a Lie algebra $L$ is {\it
perfect\,} if $[L,L]=L$). Thus using a result in [8], $\LIE$ is a
{\it complete Lie algebra}, i.e., a centerless Lie algebra without
outer derivations (see also [1,\,2]).
\par
A module $V$ over \LIE\ is called a {\it weight module} if $V$
admits a finite-dimensional weight space decomposition with
respect to the Cartan subalgebra $H=\C\ptl_1\oplus\C\ptl_2$ of
\LIE, \vs{-7pt}i.e.,
$$%
V=\rb{-5pt}{\mbox{$^{\dis\ \,\oplus}_{\l\in H^*}$}}V_\l,\ \ \
V_\l=\{v\in V\,|\,hv=\l(h)v\},\ \ {\rm dim\,}V_\l<\infty,\ \
\forall\,\l\in H^*,\vs{-10pt}\eqno(1.6)
$$%
where $H^*$ is the dual space of $H$. A weight module $V$ is
called
%a {\it uniformly bounded module} if there exists $K>0$ such
%that ${\rm dim\,}V_\l\le K$ for all $\l\in H^*$;
a {\it module of
the intermediate series} if ${\rm dim\,}V_\l\le1$ for all $\l\in
H^*$.
\par
The classification of representations, esp., irreducible
representations, of $\LIE$ is definitely an important problem and
interesting as well, not only because $\C_q$ plays an important
role in quantum groups, but also it provides examples of
irreducible representations of some complete Lie algebras which
are not simple Lie algebras. Zhang and Zhao in [9] in fact
constructed a class of modules of the intermediate series over
$\LIE$ when $q$ is generic (i.e., not a root of unity). In this
paper, we shall discuss modules of the intermediate series over
\LIE\ when $q$ is a root of unity. In particular, we classify the
modules of the intermediate series over
the derivation algebra \LIE\ of the quantum plane
for the special case when $q=-1$, i.e.,
the {\it localization of the quantum plane at $q=-1$}.
Our purpose is to understand a little bit about the
modules of the intermediate series over \LIE\ when $q$ is a root
of unity.
\par
In Section 2, after collecting some necessary information, we give
the proof of our main result (Theorem 2.3).
\par\ \vs{5pt}\par\cl{\bf 2. \ NOTATIONS AND MAIN RESULT}
\par%
Suppose $q$ is an $N$-th primitive root of unity for some $N>1$.
Then the center $Z(\C_q)=\{x^{\mm}\,|\,\mm\in N\Z^2\}$. Thus
$$%
W=Z(\C_q)\ptl_1\oplus Z(\C_q)\ptl_2={\rm
span}\{x^{\mm}\ptl_s\,|\,\mm\in N\Z^2,\,s=1,2\},\eqno(2.1)
$$%
is a rank 2 (classical) Witt algebra. We identify ${\rm
ad\,}\C'_q$ with $\C'_q$. Then by (1.4),
$$%
B=\{x^{\mm},x^{\nn}\ptl_s\,|\,\mm\in\Z^2\bs N\Z^2,\,
\nn\in N\Z^2,\,s=1,2\},$$%
is a basis of \LIE.\par%
Suppose $V$ as in (1.6) is an indecomposable weight module over
\LIE. We identify $H^*$ with $\C^2$ by the correspondence
$\l\rla(\l_1,\l_2)$, where $\l_s=\l(\ptl_s)$ for $\l\in
H^*,s=1,2$. For $\a\in\C^2$, denote
$V(\a)=\oplus_{\l\in\Z^2}V_{\a+\l}$. Clearly, $V(\a)$ is a
submodule of $V$ and $V$ is the direct sum of all different
$V(\a)$. Thus, there exists $\a\in\C^2$ such that $V=V(\a)$. For
any representative element $\ii\in\Z^2/N\Z^2$, denote
$$V[\ii]=\oplus_{\mm\in N\Z^2}V_{\a+\ii+\mm}.\eqno(2.2)$$
Obviously, $V[\ii]$ is a $W$-submodule
of $V$. \par%
Now suppose $V$ is a module of the intermediate series. Then each
$V[\ii]$ is a $W$-module of the intermediate series. Without loss
of generality, we can assume that each $V[\ii]$ is a nontrivial
$W$-module (otherwise from the discussion below, we shall see that
$V$ is a \mbox{sub-quotient} module of such a module). Fix
$\ii\in\Z^2/N\Z^2$. Then by [7], each nonzero weight of $V[\ii]$
has the
same multiplicity 1, i.e., %
$${\rm dim\,}V_{\a+\ii+\mm}=1\mbox{ \ \ for \
all \ \ }\mm\in N\Z^2\mbox{ \ \ with \ \ }\a+\ii+\mm\ne0.\eqno(2.3)$$%
Note that for any $a,b\in\C$ such that $a,b$ are $\Q$-linearly
independent, the subspace $W_1={\rm
span}\{x^{\mm}(a\ptl_1+b\ptl_2)\,|\,\mm\in N\Z^2\}$ of $W$ is a
rank 1 (generalized) Witt algebra (a centerless higher rank
Virasoro algebra). From this, one can either prove as in [5,\,6],
or directly use a result in [11] to obtain that there exists some
$b_{\ii}\in\C$, and one can choose a basis $v_{\ii+\mm}$ of
$V_{\a+\ii+\mm}$ for $\mm\in N\Z^2$ with $\a+\ii+\mm\ne0$ such
that
$$%
(x^\mm\ptl) v_{\ii+\nn}=\la\ptl,\a+\ii+\nn+b_{\ii\,}\mm\ra
v_{\ii+\mm+\nn}\mbox{ \ for \
}\a+\ii+\nn\ne0\ne\a+\ii+\mm+\nn,\eqno(2.4)
$$%
where $\mm,\nn\in N\Z^2$, and where $\la\cdot,\cdot\ra$ is the
bilinear form on $H\times\C^2$ defined by
$$%
\la\ptl,\l\ra=a_1\l_1+a_2\l_2\mbox{ \ \ for \ \
}\ptl=a_1\ptl_1+a_2\ptl_2\in
H,\,\,\l=(\l_1,\l_2)\in\C^2.\eqno(2.5)
$$%
Note that the representative element $\ii\in\Z/N\Z^2$ does not
need to be in the distinguished set $\{\mm\in\Z^2\,|\,0\le
m_1,m_2<N\}$ since if $\ii$ is replaced by $\ii+\mm$ for any
$\mm\in N\Z^2$, (2.4) still holds (where $b_\ii$ only depends on
the set $\Z^2/N\Z^2$). By choosing all $\ii\in\Z^2/N\Z^2$, the
above in fact defines $v_\mm$ for all $\mm\in\Z^2$ with
$\a+\mm\ne0$. In the following discussion, we shall use the {\bf
convention} that when a notation $v_\mm$ appears in an expression
with the index $\mm\in\Z^2$, we always
assume that $\a+\mm\ne0$. \par%
For any $\mm\in\Z^2\bs N\Z^2$, clearly $x^\mm
V_{\a+\ii+\nn}\subset V_{\a+\ii+\mm+\nn}$ since
$[\ptl_s,x^\mm]=m_sx^\mm$ for $s=1,2$, thus
$$%
x^\mm v_{\ii+\nn}=c_{\mm,\ii+\nn}v_{\ii+\mm+\nn}\mbox{ \ for \
}\nn\in N\Z^2\mbox{ \ and for some \
}c_{\mm,\ii+\nn}\in\C.\eqno(2.6)
$$%
Note that when $c_{\mm,\ii+\nn}$ appears in an expression, we
always assume that $\a+\ii+\nn\ne0\ne\a+\ii+\mm+\nn$. Let $\ptl\in
H,\,\kk\in N\Z^2$. Applying $\la\ptl,\mm\ra
x^{\kk+\mm}=[x^\kk\ptl,x^\mm]$ to $v_{\ii+\nn}$, by (2.4) and
(2.6), we obtain
$$%
\la\ptl,\mm\ra
c_{\kk+\mm,\ii+\nn}=\la\ptl,\a+\ii+\mm+\nn+b_{\ii+\mm\,}\kk\ra
c_{\mm,\ii+\nn}-
c_{\mm,\ii+\kk+\nn}\la\ptl,\a+\ii+\nn+b_{\ii\,}\kk\ra.\eqno(2.7)
$$%
Note that in (2.7), by our convention we always assume
$$\a+\ii+\nn,\ \a+\ii+\mm+\nn,\ \a+\ii+\kk+\nn,\ \a+\ii+\kk+\mm+\nn\ \ne\ 0.
\eqno(2.8)$$
{\bf Lemma 2.1}. {\it (i) $c_{\mm,\ii}=c_{\mm+\kk,\ii+\nn}$ for
all $\kk,\nn\in N\Z^2$; (ii) $b_\ii=b_{\ii+\mm}$ if
$c_{\mm,\ii}\ne0$.
}\par\ni%
{\bf Proof.} %
Suppose $\mm\in\Z^2\bs N\Z^2,\ii\in\Z^2/N\Z^2$ are fixed. We
consider the
following cases. \par%
{\it Case 1:} $b_\ii=b_{\ii+\mm}$.\vs{0pt}\par %
Denote $b=b_\ii$. Take $\ptl'=m_2\ptl_1-m_1\ptl_2\in H$. Then
$\la\ptl',\mm\ra=0$ by (2.5). For any $\l\in\C^2$, for simplicity
we denote $\l_{\ptl'}=\la\ptl,\l\ra\in\C$ and denote
$\ii_{\ptl'}^\a=\a_{\ptl'}+\ii_{\ptl'}$. Then (2.7) gives
$$
(\ii^\a_{\ptl'}+\nn_{\ptl'}+b\kk_{\ptl'})(c_{\mm,\ii+\nn}-
c_{\mm,\ii+\kk+\nn})=0.\eqno(2.9)
$$
We claim
$$c_{\mm,\ii+\nn}=c_{\mm,\ii+\kk+\nn}\mbox{ \ for \ }\kk,\nn\in
N\Z^2.\eqno(2.10)$$%
First suppose $b=1$. Choose $\wt{\kk}\in N\Z^2$ such that
$\ii^\a_{\ptl'}+\wt{\kk}_{\ptl'}\ne0$ (note that since
$\ptl'\ne0$, $\kk_{\ptl'}$ can take infinite many different values
for $\kk\in N\Z^2$). Now for any $\nn\in N\Z^2$, we take
$\kk=\wt{\kk}-\nn$ in (2.9), then we obtain
$c_{\mm,\ii+\nn}=c_{\mm,\ii+\wt{\kk}}$. Thus we have (2.10).
Similarly (2.10) holds if $b=0$.\par %
Next suppose $b\ne0,1$. Then for any $\kk,\nn\in N\Z^2$, we can
always choose some $\jj\in N\Z^2$ such that
$$%
\ii^\a_{\ptl'}+\nn_{\ptl'}+b(\jj_{\ptl'}-\nn_{\ptl'})\ne0\ne
\ii^\a_{\ptl'}+\jj_{\ptl'}+b(\kk_{\ptl'}+\nn_{\ptl'}-\jj_{\ptl'}).
$$%
Thus (2.9) gives that
$c_{\mm,\ii+\nn}=c_{\mm,\ii+\jj}=c_{\mm,\ii+\kk+\nn}$, and (2.10) holds.\par %
Now choose $\ptl\in H$ such that $\mm_\ptl\ne0$. Then (2.10) and
(2.7) show that $c_{\kk+\mm,\ii+\nn}=c_{\mm,\ii}$. This proves the
lemma in this case.\par%
 {\it Case 2:} $b_{\ii}\ne b_{\ii+\mm}$ and
$c_{\mm,\ii+\nn}=0$ for all $\nn\in N\Z^2$.
\par%
In this case, in particular (2.10) holds, which together with
(2.7) gives (i) as in case 1. For (ii), there is nothing to prove in this case.\par%
{\it Case 3:} $b_{\ii}\ne b_{\ii+\mm}$ and
$c_{\mm,\ii+\wt{\nn}}\ne0$ for some $\wt{\nn}\in N\Z^2$.
\par%
As in Case 1, take $\ptl'=m_2\ptl_1-m_1\ptl_2\in H$. Then (2.7)
gives
$$
(\ii^\a_{\ptl'}+\nn_{\ptl'}+b_{\ii+\mm\,}\kk_{\ptl'})
c_{\mm,\ii+\nn}-
(\ii^\a_{\ptl'}+\nn_{\ptl'}+b_{\ii\,}\kk_{\ptl'})c_{\mm,\ii+\kk+\nn}=0.\eqno(2.11)
$$
Denote the left-hand side of (2.11) by $d(\kk,\nn)$. Then for any
$\jj,\kk\in\Z^2$, we have a system of 3 linear equations in the 3
unknown variables
$c_{\mm,\ii+\wt{\nn}},c_{\mm,\ii+\kk+\wt{\nn}},c_{\mm,\ii+\jj+\wt{\nn}}$:$$d(\kk,\wt{\nn})=d(\jj,\wt{\nn})=d(\jj-\kk,\kk+\wt{\nn})=0.\eqno(2.12)$$
Since (2.12) has a solution $c_{\mm,\ii+\wt{\nn}}\ne0$, we obtain
that the determinant of these 3 linear equations is zero, i.e., %
\par\ni\hs{2pt}$ %
\begin{array}{ll}
&\left|\!\begin{array}{ccc}
\ii^\a_{\ptl'}\!+\!\wt{\nn}_{\ptl'}\!+\!b_{\ii+\mm\,}\kk_{\ptl'}&
-\ii^\a_{\ptl'}\!-\!\wt{\nn}_{\ptl'}\!-\!b_{\ii\,}\kk_{\ptl'}&0\vs{4pt}\\\ii^\a_{\ptl'}\!+\!\wt{\nn}_{\ptl'}\!+\!b_{\ii+\mm\,}\jj_{\ptl'}&0&
-\ii^\a_{\ptl'}\!-\!\wt{\nn}_{\ptl'}\!-\!b_{\ii\,}\jj_{\ptl'}\vs{4pt}\\
0&\!\ii^\a_{\ptl'}\!+\!\kk_{\ptl'}\!+\!\wt{\nn}_{\ptl'}\!+\!b_{\ii+\mm}
(\jj_{\ptl'}\!-\!\kk_{\ptl'})&
-\ii^\a_{\ptl'}\!-\!\kk_{\ptl'}\!-\!\wt{\nn}_{\ptl'}\!-\!b_{\ii}
(\jj_{\ptl'}\!-\!\kk_{\ptl'})
\end{array}\!\!\right|\vs{12pt}\\
=\!\!\!&(\jj_{\ptl'}-\kk_{\ptl'})\kk_{\ptl'} (b_{\ii+\mm}-b_{\ii})
   ((\ii^\a_{\ptl'}+\wt{\nn}_{\ptl'})
   (b_{\ii}+b_{\ii+\mm}-1)+\jj_{\ptl'}b_{\ii}b_{\ii+\mm})\vs{6pt}\\
=\!\!\!&
0,\end{array} %
$\hfill(2.13)\par\ni %
for all $\jj,\kk\in\Z^2$. Since $\jj_{\ptl'},\kk_{\ptl'}$ can take
infinite many different values, (2.13) shows that $b_\ii
b_{\ii+\mm}=0$. Note that in (2.12), if we replace
$\wt{\nn},\jj,\kk$ by $\wt{\nn}+\jj,-\jj,\kk-\jj$ respectively, we
have another system of 3 linear equations in the 3 unknown
variables
$c_{\mm,\ii+\wt{\nn}},c_{\mm,\ii+\kk+\wt{\nn}},c_{\mm,\ii+\jj+\wt{\nn}}$,from
which, as discussion in (2.13), we obtain that
$b_{\ii}+b_{\ii+\mm}=1$. Thus we obtain
$$%
b_{\ii}=0,\,b_{\ii+\mm}=1,\mbox{ \ \ or \ \
}b_{\ii}=1,\,b_{\ii+\mm}=0. \eqno(2.14)$$ %
\par%
If necessary by replacing $\ii,\mm$ by $\ii+\mm,-\mm$
respectively, we can suppose $b_\ii=0,\,b_{\ii+\mm}=1$.\par%
If $\ii^\a_{\ptl'}+\wt{\nn}_{\ptl'}=0$, then by choosing $\kk\in
N\Z^2$ with $\kk_{\ptl'}\ne0$, from (2.11), we obtain
$c_{\mm,\ii+\wt{\nn}}=0$, a contradiction. Thus
$\ii^\a_{\ptl'}+\wt{\nn}_{\ptl'}\ne0$, and (2.11) gives
$$
c_{\mm,\ii+\kk}=(\ii^\a_{\ptl'}+\wt{\nn}_{\ptl'})^{-1}
(\ii^\a_{\ptl'}+\kk_{\ptl'})c_{\mm,\ii+\wt{\nn}_{\ptl'}}\mbox{ \
for \ }\kk\in N\Z^2.\eqno(2.15)
$$
Choose $\ptl''\in H$ such that $\mm_{\ptl''}=1$. Then (2.7) and
(2.15) give
$$
c_{\kk+\mm,\ii+\nn}=((\ii^\a_{\ptl''}+1+\nn_{\ptl''}+\kk_{\ptl''})(\ii^\a_{\ptl'}
+\nn_{\ptl'})-(\ii^\a_{\ptl''}+\nn_{\ptl''})(\ii^\a_{\ptl'}+\kk_{\ptl'}+\nn_{\ptl'}))
(\ii^\a_{\ptl'}+\wt{\nn}_{\ptl'})^{-1}c_{\mm,\ii+\wt{\nn}}.
\eqno(2.16)
$$
In (2.7), replacing $\ptl,\mm,\kk$ by $\ptl'',\kk+\mm,-\kk$
respectively, noting that $b_{\ii+\kk+\mm}=b_{\ii+\mm}$ for
$\kk\in N\Z^2$, we obtain
$$%
(\kk_{\ptl''}+1)c_{\mm,\ii+\nn}=(\ii^\a_{\ptl''}+1+\nn_{\ptl''})c_{\kk+\mm,\ii+\nn}-
(\ii^\a_{\ptl''}+\nn_{\ptl''})c_{\kk+\mm,\ii-\kk+\nn}.\eqno(2.17)
$$
Using (2.15) and (2.16) in (2.17), we obtain
$$%
(\kk_{\ptl'}-\kk_{\ptl''})(\ii^\a_{\ptl'}+\nn_{\ptl'})
(\ii^\a_{\ptl'}+\wt{\nn}_{\ptl'})^{-1}c_{\mm,\ii+\wt{\nn}}=0.\eqno(2.18)$$
Clearly we can choose $\kk\in N\Z^2$ such that
$\kk_{\ptl'}\ne\kk_{\ptl''}$ since $\ptl'\ne\ptl''$. Then (2.18)
shows that $c_{\mm,\ii+\wt{\nn}}=0$, a contradiction.
Thus this case does not occur, and the lemma is proved.\qed\par%
Now for $\ii\in\Z^2/N\Z^2$, denote %
$$\wt V[\ii]={\rm span}\{v\in
V[\jj]\,|\,\jj\in\Z^2/N\Z^2,\,b_\jj=b_\ii\}.\eqno(2.19)$$%
Then (2.2) and Lemma 2.1(ii) show that $\wt V[\ii]$ is a
(\LIE)-submodule of $V$, and obviously, $V$ is the direct sum of
all different $\wt V[\ii]$. Since $V$ is indecomposable, we obtain
$V=\wt V[\ii]$ for some $\ii\in\Z^2/N\Z^2$, i.e., $b_\ii=b$ is a
constant for $\ii\in\Z^2/N\Z^2$.
\par
Set $c_{\mm,\ii}=0$ if $\mm\in N\Z^2,\,\ii\in\Z^2$. Lemma 2.1(i)
shows that $c_{\mm,\ii}$ is a function on
$(\Z^2/N\Z^2)\times(\Z^2/N\Z^2)$. Applying (1.1) to $v_{\ii}$, we
obtain a system of equations in unknown variables
$c_{\mm,\ii},\,\mm,\ii\in\Z^2/N\Z^2$:
$$%
c_{0,\ii}=0,\,\,\
c_{\mm,\nn+\ii}c_{\nn,\ii}-c_{\nn,\mm+\ii}c_{\mm,\ii}=
(q^{m_2n_1}-q^{m_1n_2})c_{\mm+\nn,\ii}\mbox{ \ \ for \ \
}\mm,\nn,\ii\in\Z^2/N\Z^2.\eqno(2.20)
$$%
From the discussion above, we see that the classification of
indecomposable (\LIE)-module $V$ of the intermediate series is
equivalent to solving the system (2.20).
\par%
For $a=0,1,\,\a\in\C^2,\,b\in\C$, we define a (\LIE)-module
$A_{a,\a,b}$ of the intermediate series as follows: It has a basis
$\{v_\kk\,|\,\kk\in\Z^2\}$ such that
$$%
x^\mm v_\kk=(q^{m_2k_1}-a q^{m_1k_2})v_{\mm+\kk},\ \ \
(x^\nn\ptl)v_\kk=\la\ptl,\a+\kk+b\nn\ra v_{\nn+\kk},\eqno(2.21)
$$%
for $\mm\in\Z^2\bs N\Z^2,\,\nn\in N\Z^2,\,\kk\in\Z^2,\,\ptl\in H$.
\par\ni%
{\bf Proposition 2.2}. {\it (i) $A_{0,\a,b}$ is an irreducible
$(\LIE)$-module; (ii) $A_{1,\a,b}=A'_{1,\a,b}\oplus A''_{1,\a,b}$
is a direct sum of two $(\LIE)$-submodules, where
$A'_{1,\a,b}={\rm span}\{v_\kk\,|\,\kk\in\Z^2\bs N\Z^2\}$ is
irreducible and $A''_{1,\a,b}={\rm span}\{v_\kk\,|\,\kk\in
N\Z^2\}$ is a $(\LIE)$-module such that the action of ${\rm
ad\,}\C'_q$ is trivial (thus $A''_{1,\a,b}$ is simply a
$W$-module, which is irreducible if and only if $\a\notin N\Z^2$
or $b\ne0,1$).
\par\ni%
Proof.} %
It is straightforward.\qed\par\ni%
{\bf Theorem 2.3}. {\it Suppose $q=-1$ (i.e., $N=2$). An
irreducible $(\LIE)$-module $V$ of the intermediate series is a
sub-quotient module of
$A_{a,\a,b}$ for some $a=0,1,\,\a\in\C^2,\,b\in\C$.\par\ni%
Proof.} %
We can suppose $V$ is a {\it faithful module} (i.e., every nonzero
element of \LIE\ acts nontrivially on $V$), otherwise the only
proper ideal ${\rm ad\,}\C'_q$ of \LIE\ must act trivially on $V$
and then $V$ is simply a
$W$-module, and so the result follows from [11].\par%
{\it Case 1}: Suppose each $V[\ii]$ is a nontrivial $W$-module. \par%
We want to prove that by a suitable choice of basis elements
$v_\kk,\,\kk\in\Z^2$ satisfying (2.6), we have
$$
c_{\mm,\kk}=(-1)^{m_2k_1}\mbox{ \ for \
}\mm\in\Z^2\bs2\Z^2,\,\kk\in\Z^2.\eqno(2.22)
$$
In fact, by (2.20), it suffices to prove (2.22) for
$\mm=(1,0),(0,1)$, i.e., we only need to solve 8 unknown variables
$$c_{\mm,\kk}\mbox{ \ \ for \ \ }\mm=(1,0),(0,1),\,\kk\in\Z^2/2\Z^2.\eqno(2.23)$$
First, we have $x_1V[(0,0)]\ne0$ or $x_2V[(0,0)]\ne0$, otherwise
$V[(0,0)]$ is a proper submodule. \par%
{\it Subcase (a)}: Suppose $x_1V[(0,0)]\ne0$ and $x_2V[(0,0)]\ne0$. \vs{0pt}\par%
We have $x_1V[(0,1)]\ne0$ or $x_2V[(1,0)]\ne0$, otherwise
$V[(0,0)]\oplus V[(1,0)]\oplus V[(0,1)]$ is a proper submodule. By
symmetry, we can suppose $x_1V[(0,1)]\ne0$. Thus, by rescaling
basis elements, we can suppose that (2.22) holds for
$\mm=(1,0),\,\kk=(0,0),(0,1)$ and $\mm=(0,1),\,\kk=(1,0)$.
\par%
Thus 3 of 8 variables in (2.23) are known, using this, one can
solve the system (2.20) to obtain that (2.22) holds in general
(this can be easily done by the computer software {\it
Mathematika}).
\par{\it Subcase (b)}: By symmetry, we can suppose $x_1V[(0,0)]\ne0$ and
$x_2V[(0,0)]=0$. \par%
We have $x_2V[(1,0)]\ne0$ and $x_1V[(1,1)]\ne0$, otherwise either
$V[(0,0)]\oplus V[(1,0)]$ or $V[(0,0)]\oplus V[(1,0)]\oplus
V[(1,1)]$ is a proper submodule. Thus similarly, we can suppose
that (2.22) holds for $\mm=(1,0),\,\kk=(0,0),(1,1)$ and
$\mm=(0,1),\,\kk=(1,0)$, and similarly, using this, one can solve
the system (2.20) to
obtain that (2.22) holds in general.\par %
Thus (2.22) holds and $V$ is the module $A_{0,\a,b}$ (we remark
that in case $\a\in\Z^2$, the above discussion only shows that we
have (2.21) (with $a=0,q=-1,N=2$) for all $\mm,\nn,\kk$ such that
$v_{-\a}$ is not involved in the expression; then we can use the
first equation of (2.21) to define $v_{-\a}$ by
$v_{-\a}=x^{(1,0)}v_{-\a-(1,0)}$ and show that (2.21) holds in
general). This proves the theorem in this case.\par%
{\it Case 2}:  $V[\ii]$ is a trivial $W$-module for some
$\ii\in\Z^2/2\Z^2$. \par%
By replacing $\a$ by $\a+\ii$, we can suppose $\ii=(0,0).$ We
claim that $V[\jj]$ is nontrivial for $\jj=(1,0),(0,1),(1,1)$: If
$V[(1,0)]$ is trivial, then we can check that $V'=V[(0,1)]\oplus
V[(1,1)]$ is a $(\LIE)$-submodule and $x_2$ acts trivially on
$V'$; if $V[(1,1)]$ is trivial, then $V'=V[(1,0)]\oplus V[(0,1)]$
is a submodule and $x_1,x_2$ act trivially on $V'$; if $V[(0,1)]$
is trivial, then $V'=V[(1,0)]\oplus V[(1,1)]$ is a submodule and
$x_1$ acts trivially on $V'$; all lead to a contradiction with the
assumption that $V$ is irreducible and faithful.\par%
We can prove similarly that $x_1V[(1,1)]\ne0$ and
$x_2V[(1,1)]\ne0$. Thus by rescaling basis elements, we can
suppose
$$%
c_{\mm,\kk}=(-1)^{m_2k_1}-(-1)^{m_1k_2},
\eqno(2.24)$$%
holds for $\mm=(1,0),\,(0,1),\,\kk=(1,1)$. Since
$x_1V[(1,0)]=x_2V[(0,1)]=x_1V[(0,0)]=x_2V[(0,0)]=0$, we can
suppose $c_{\mm,\kk}=0$ for $\mm=\kk=(1,0)$, $\mm=\kk=(0,1)$ and
$\mm=(1,0),\,(0,1),\,\kk=(0,0)$, i.e., (2.24) also holds for
$\mm=\kk=(1,0)$, $\mm=\kk=(0,1)$ and
$\mm=(1,0),\,(0,1),\,\kk=(0,0)$. Thus 6 of 8 variables in (2.23)
are known, using this, one can solve the system (2.20) to obtain
that (2.24) holds in general (again, this can be easily done by
the computer software {\it Mathematika}). Thus $V$ is the module
$A'_{1,\a,b}$. This completes the proof of the theorem.\qed\par
\vskip 13pt \small \par\ni {\bf Acknowledgements}
\par\ni
This work is supported by NSF grant no.~10171064, 10471091 and no.~10271076
of China, two grants ``Excellent Young Teacher Program'' and
``Trans-Century Training Programme Foundation for the Talents''
from Ministry of Education of China and the Foundation of Jiangsu
Educational Committee. They are grateful to Professor K. Zhao for
many useful discussions on the subject.
 \vskip 13pt
\par\ni
{\bf References}%
\vs{0pt}\par\ni\hi3.5ex\ha1 %
[1]~Jiang, C.~and Meng,~D.: The derivation algebra of the
associative algebra $\C_q[X,Y,X^{-1},Y^{-1}]$, {\it Comm.~Alg.}
{\bf26} (1998), 1723-1736.%
  \par\ni\hi3.5ex\ha1 %
[2]~Jiang, C.~and Meng,~D.: The derivation algebra and the
universal central extension of the $q$-analog of the Virasoro-like
algebra, {\it Comm.~Alg.} {\bf26} (1998), 1335-1346.%
  \par\ni\hi3.5ex\ha1 %
[3]~Kassel, C.: {\it Quantum Groups,} Springer-Verlag, 2000.%
  \par\ni\hi3.5ex\ha1 %
[4]~Kirkman, E., Procesi, C.~and Small, L.: A $q$-analog for the
Virasoro algebra, {\it Comm.~Alg.} {\bf22} (1994), 3755-3774.%
  \par\ni\hi3.5ex\ha1 %
[5]~Su,~Y.: Harish-Chandra modules of the intermediate series over
the high
 rank Virasoro algebras and high rank super-Virasoro algebras,
 {\it J. Math. Phys.} {\bf35} (1994), 2013-2023.
 \par\ni\hi3.5ex\ha1 %
[6]~Su,~Y.: Simple modules over the high rank Virasoro algebras,
{\it Comm.~Alg.} {\bf29} (2001), 2067-2080.
  \par\ni\hi3.5ex\ha1 %
[7]~Su,~Y.: Classification of Harish-Chandra modules over the
higher rank Virasoro algebras, {\it Comm.~Math.~Phys.}
{\bf240} (2003), 539-551.
  \par\ni\hi3.5ex\ha1 %
[8]~Su,~Y.~and Zhu,~L.: Derivation algebras of centerless perfect
Lie algebras are complete, {\it J.~Alg.}, in press.
  \par\ni\hi3.5ex\ha1 %
[9]~Zhang,~H.~and Zhao,~K.: Representations of the Virasoro-like
Lie algebra  and its $q$-analog, {\it Comm.~Alg.} {\bf24} (1996),
4361-4372.%
  \par\ni\hi3.5ex\ha1 %
[10]~Zhao,~K.: The $q$-Virasoro-like algebra,  {\it J.~Alg.}
{\bf188} (1997), 506-512.%
  \par\ni\hi3.5ex\ha1 %
[11]~Zhao,~K.: Harish-Chandra modules over generalized Witt
algebras, to appear.
\end{document}